\documentclass[12pt]{amsart}
\usepackage{amssymb}
\usepackage{amsfonts}
\usepackage{color}
\newtheorem{thm}{Theorem}[section]

\theoremstyle{definition}

\theoremstyle{remark}

\newcommand{\X}{\mathfrak X}
\newcommand{\B}{\mathfrak B}

\newcommand{\bb}{\mathfrak b}
\begin{document}

\title[K\"AHLER MANIFOLDS ADMITTING A FLAT COMPLEX...]
{K\"AHLER MANIFOLDS ADMITTING A FLAT COMPLEX CONFORMAL CONNECTION}%
\author{G. Ganchev and V. Mihova}%
\address{Bulgarian Academy of Sciences, Institute of Mathematics and
Informatics, Acad. G. Bonchev Str. bl. 8, 1113 Sofia, Bulgaria}%
\email{ganchev@math.bas.bg}%
\address{Faculty of Mathematics and Informatics, University of Sofia,
J. Bouchier Str. 5, (1164) Sofia, Bulgaria}
\email{mihova@fmi.uni-sofia.bg}
\subjclass{Primary 53B35, Secondary 53C25}%
\keywords{Complex conformal connection, Bochner-K\"ahler manifolds with
special scalar distribution, warped product K\"ahler manifolds.}%

\begin{abstract}
\vskip 2mm
We prove that any K\"ahler manifold admitting a flat complex conformal
connection is a Bochner-K\"ahler manifold with special scalar
distribution and zero geometric constants. Applying the local structural
theorem for such manifolds we obtain a complete description of
the K\"ahler manifolds under consideration.
\end{abstract}
\maketitle
\section{Introduction}
Let $(M,g,J)\;(\dim \,M=2n)$ be a K\"ahler manifold with complex structure
$J$, metric $g$, Levi-Civita connection $\nabla$, curvature tensor $R$,
Ricci tensor $\rho$ and scalar curvature $\tau$.
The Bochner curvature tensor $B(R)$ is given by
$$\begin{array}{l}
B(R)(X,Y)Z=R(X,Y)Z\\
[2mm]
-Q(Y,Z)X+Q(X,Z)Y-g(Y,Z)Q(X)+g(X,Z)Q(Y)\\
[2mm]
-Q(JY,Z)JX+Q(JX,Z)JY+2Q(JX,Y)JZ\\
[2mm]
-g(JY,Z)JQ(X)+g(JX,Z)JQ(Y)+2g(JX,Y)JQ(Z), \quad X,Y,Z\in {\X}M,
\end{array} $$
where $Q(X,Y)=\displaystyle{\frac{1}{2(n+2)}\,\rho(X,Y)
-\frac{\tau}{8(n+1)(n+2)}\,g(X,Y)}$ and $Q(X)$ is the corresponding
tensor of type (1.1).

The manifold is said to be {\it Bochner flat (Bochner-K\"ahler)} if its
Bochner curvature tensor vanishes identically, i.e.
$$\begin{array}{l}
R(X,Y)Z=\\
[2mm]
Q(Y,Z)X-Q(X,Z)Y+g(Y,Z)Q(X)-g(X,Z)Q(Y)\\
[2mm]
+Q(JY,Z)JX-Q(JX,Z)JY-2Q(JX,Y)JZ\\
[2mm]
+g(JY,Z)JQ(X)-g(JX,Z)JQ(Y)-2g(JX,Y)JQ(Z), \quad X,Y,Z\in {\X}M,
\end{array} \leqno{(1.1)}$$

For any real $\mathcal{C}^{\infty}$ function $u$ on $M$
we denote $\omega=du$ and $P= grad\,u$.

In \cite{Y} Yano introduced on a K\"ahler manifold a complex conformal
connection and proved
\vskip 2mm
{\bf Theorem A.} {\it If in a $2n$-dimensional\, $(n\geq 2)$ \,K\"ahler
manifold there exists a scalar function $u$ such that the
complex conformal connection
$$\begin{array}{ll}
\mathcal{D}_XY=&\nabla_XY+ \omega(X)Y+\omega(Y)X-g(X,Y)P\\
[2mm]
&-\omega(JX)JY-\omega(JY)JX-g(JX,Y)JP, \quad X,Y\in{\X}M,
\end{array}$$
is of zero curvature, then the
Bochner curvature tensor of the manifold vanishes.}
\vskip 2mm
In \cite{S} Seino proved the inverse
\vskip 2mm
{\bf Theorem B.} {\it In a K\"ahlerian space with vanishing Bochner
curvature tensor if there exists a non-constant function $u$ satisfying
the equality
$$(\nabla_X\omega)(Y)+2\omega(JX)\omega(JY)+\omega(P)g(X,Y)=0,$$
then the complex conformal connection is of zero curvature.}
\vskip 2mm
In this paper we prove
\vskip 2mm
{\bf Theorem 3.1.} {\it A K\"ahler manifold $(M,g,J) \;
(\dim M = 2n \geq 6)$ admits a flat complex conformal connection
if and only if it is a Bochner-K\"ahler manifold whose scalar distribution
$D_{\tau}$ is a $B_0$-distribution with function $a+k^2=0$ and
geometric constants ${\B}=\frak b_0=0.$}
\vskip 2mm
Applying the local structural theorem \cite{GM2} for Bochner-K\"ahler
manifolds whose scalar distribution is a $B_0$-distribution, we
describe locally all K\"ahler manifolds admitting a flat complex
conformal connection.

\section{Preliminaries}
Let $(M,g,J) \, (\dim M=2n)$ be a K\"ahler manifold with metric $g$,
complex structure $J$ and Levi-Civita connection $\nabla$. We denote by
${\X}M$ the Lie algebra of all $\mathcal{C}^{\infty}$ vector fields on $M$.
The fundamental K\"ahler form $\Omega$ is defined as follows
$$\Omega(X,Y)=g(JX,Y), \quad X,Y \in {\X}M.$$

For any $\mathcal{C}^{\infty}$ real function $u$ on $M$ we consider
the conformal metric $\bar g= e^{2u}g$. We denote the 1-form
$\omega:= du$ and $P:= grad\, u$ with respect to the metric $g$.
Then $(M,\bar g,J)$ is a locally conformal K\"ahler manifold, or a
$W_4$-manifold in the classification scheme of  \cite{GH}. The
fundamental K\"ahler form and the Lee form of the structure
$(\bar g,J)$ are $\bar \Omega (X,Y)=\bar g(JX,Y),
\; X,Y \in {\X}M$ and $\bar \omega =2\omega = 2du$, respectively.
The Lee vector $\bar P$ corresponding to $\bar \omega$ with respect to
the metric $\bar g$ is $\bar P= 2e^{-2u}P$.

The unique linear connection $\mathcal{D}$ with torsion $\mathcal{T}$
satisfying the conditions:
$$\begin{array}{l}
1)\, \mathcal{D}J=0;\\
[2mm]
2)\, \mathcal{D}\bar g=0;\\
[2mm]
3)\, \mathcal{T}=-\,\bar\Omega \otimes J\bar P\end{array} \leqno{(2.1)}$$
is said to be a {\it complex conformal connection}  \cite{Y}.

In terms of the K\"ahler structure $(g,J)$ $\mathcal{D}$ is given by
$$\begin{array}{ll}
\mathcal{D}_XY=&\nabla_XY+ \omega(X)Y+\omega(Y)X-g(X,Y)P\\
[2mm]
&-\omega(JX)JY-\omega(JY)JX-g(JX,Y)JP, \quad X,Y\in{\X}M.
\end{array}\leqno{(2.2)}$$
The conditions (2.1) in terms of the K\"ahler structure $(g,J)$
become
$$\begin{array}{l}
1)\, \mathcal{D}J=0;\\
[2mm]
2)\, \mathcal{D}g=-2\,\omega \otimes g;\\
[2mm]
3)\, \mathcal{T}=-2\,\Omega \otimes JP.\end{array} \leqno{(2.3)}$$

Denote by $\mathcal{R}$ the curvature tensor of the complex conformal
connection $\mathcal {D}$. Taking into account (2.2) we have the
relation between $R$ and $\mathcal{R}$:

$$\begin{array}{l}
\mathcal{R}(X,Y)Z=R(X,Y)Z\\
[1mm]
-\{(\nabla_Y\omega)(Z)-\omega(Y)\omega(Z)+\omega(JY)\omega(JZ)
+\displaystyle{\frac{1}{2}}\,\omega(P)g(Y,Z)\}X\\
[2mm]
+\{(\nabla_X\omega)(Z)-\omega(X)\omega(Z)+\omega(JX)\omega(JZ)
+\displaystyle{\frac{1}{2}}\,\omega(P)g(X,Z)\}Y\\
[2mm]
-g(Y,Z)\{\nabla_XP-\omega(X)P-\omega(JX)JP+\displaystyle{\frac{1}{2}}
\omega(P)X\}\\
[2mm]
+g(X,Z)\{\nabla_YP-\omega(Y)P-\omega(JY)JP+\displaystyle{\frac{1}{2}}
\omega(P)Y\}\\
[2mm]
+\{(\nabla_Y\omega)(JZ)-\omega(Y)\omega(JZ)-\omega(JY)\omega(Z)
+\displaystyle{\frac{1}{2}}\,\omega(P)g(Y,JZ)\}JX\\
[2mm]
-\{(\nabla_X\omega)(JZ)-\omega(X)\omega(JZ)-\omega(JX)\omega(Z)
+\displaystyle{\frac{1}{2}}\,\omega(P)g(X,JZ)\}JY\\
[2mm]
+g(Y,JZ)\{\nabla_XJP-\omega(X)JP+\omega(JX)P+\displaystyle{\frac{1}{2}}
\omega(P)JX\}\\
[2mm]
-g(X,JZ)\{\nabla_YJP-\omega(Y)JP+\omega(JY)P+\displaystyle{\frac{1}{2}}
\omega(P)JY\}\\
[2mm]
-(\nabla_X\omega)(JY)JZ+(\nabla_Y\omega)(JX)JZ+
2g(X,JY)\{\omega(JZ)P+\omega(Z)JP\}
\end{array}\leqno{(2.4)}$$
for all $X,Y,Z\in{\X}M.$

From (2.4) it follows that the curvature tensor $\mathcal {R}$ satisfies
the first Bianchi identity (i.e. $\mathcal {R}$ is a K\"ahler tensor) if
and only if \cite{S}:
$$(\nabla_X\omega)(Y)+2\omega(JX)\omega(JY)+\omega(P)g(X,Y)=0,\quad
X,Y\in{\X}M,\leqno{(2.5)}$$
which is equivalent to the condition
$$\mathcal{D}_XP=0, \quad X\in{\X}M.$$

If the 1-form $\omega$ satisfies (2.5), then (2.4) becomes
$$\begin{array}{l}
\mathcal{R}(X,Y)Z=R(X,Y)Z\\
[2mm]
+L(Y,Z)X-L(X,Z)Y+g(Y,Z)L(X)-g(X,Z)L(Y)\\
[2mm]
+L(JY,Z)JX-L(JX,Z)JY-2L(JX,Y)JZ\\
[2mm]
+g(JY,Z)JL(X)-g(JX,Z)JL(Y)-2g(JX,Y)JL(Z), \quad X,Y,Z\in {\X}M,
\end{array} \leqno{(2.6)}$$
where $L(X,Y)=\omega(X)\omega(Y)+\omega(JX)\omega(JY)+
\displaystyle{\frac{1}{2}}\,\omega(P)g(X,Y)$ and $L(X)$ is the
corresponding tensor of type (1,1) with respect to the K\"ahler
metric $g$.

If $(M,g,J)$ admits a flat complex conformal connection (2.2),
then $\mathcal{R}$ satisfies the first Bianchi identity, i.e.
(2.5) holds good. Then (2.6) implies that the K\"ahler manifold is
Bochner flat.

Conversely, if $(M,g,J)$ admits a 1-form $\omega$ satisfying (2.5),
then (2.4) becomes (2.6). The condition $(M,g,J)$ is Bochner
flat implies that ${\mathcal R}=0$, i.e. the complex conformal connection
(2.2) is flat.

\section{A Curvature characterization of K\"ahler manifolds admitting
flat complex conformal connection}
For any Bochner-K\"ahler manifold $(M,g,J)$ in \cite{GM2} we proved that
$$\begin{array}{ll}
(\nabla_X\,\rho)(Y,Z)=&\displaystyle{\frac{1}{4(n+1)}\,\{2d\tau(X)g(Y,Z)
+d\tau(Y)g(X,Z)+d\tau(Z)g(X,Y)}\\
[3mm]
&+d\tau(JY)g(X,JZ)+d\tau(JZ)g(X,JY)\}, \quad X,Y,Z \in {\X}M.
\end{array}\leqno (3.1)$$
This equality shows that
the conditions $\tau = const$ and $\nabla\rho = 0$ are equivalent on a
Bochner-K\"ahler manifold. Because of the structural theorem in
\cite {TL} the case $B(R)=0, \; d\tau =0$, can be considered as
well-studied.

We consider Bochner-K\"ahler manifolds satisfying the condition
$d\tau \neq 0$ for all points $p \in M.$
This condition allows us to introduce the frame field
$$\left\{ \xi = \frac{grad \, \tau}{\Vert d\tau \Vert},
\quad J\xi = \frac{Jgrad \, \tau}{\Vert d\tau \Vert} \right\}$$
and the $J$-invariant distributions $D_{\tau}$ and
$D^{\perp}_{\tau}=span\{\xi,J\xi\}$.

Thus our approach to the local theory of Bochner-K\"ahler manifolds is
to treat them as K\"ahler manifolds $(M,g,J,D_{\tau})$ endowed with
a $J$-invariant distribution $D_{\tau}$ generated by the K\"ahler
structure $(g,J)$. We call this distribution
{\it the scalar distribution} of the manifold \cite{GM2}.

A J-invariant distribution
$D_{\tau},\, (D_{\tau}^{\perp} = span \{\xi, J\xi \})$ is said to be a
$B_0$-distribution \cite {GM1} if $\dim M=2n\geqq 6$ and
$$\begin{array}{l}
i) \displaystyle{ \quad \nabla _{x_0} \xi =
\frac{k}{2}\,x_0, \quad k \neq 0,}\quad x_0 \in D_{\tau},\\
[2mm]
ii) \quad \nabla_{J\xi}\xi=-p^*J\xi,\\
[2mm]
iii)\quad \nabla _{\xi} \xi = 0,
\end{array}$$
where $k$ and $p^*$ are functions on $M$.

The above conditions are equivalent to the equalities
$$\begin{array}{l}
\displaystyle{\nabla_X\xi=\frac{k}{2}\{X-\eta(X)\xi+\eta(JX)J\xi\}
+p^*\eta(JX)J\xi, \quad X\in {\X}M,}\\
[2mm]
\displaystyle{dk = \xi(k)\,\eta, \quad p^* = -\frac{\xi(k)+k^2}{k}.}
\end{array} \leqno (3.2)$$

In \cite{GM2} we have shown that
$${\B}=\Vert \rho \Vert^2 - \frac{\tau^2}{2(n+1)} +
\frac{\Delta \tau}{n+1}\leqno{(3.3)}$$
is a constant on any Bochner-K\"ahler manifold. We call this constant
{\it the Bochner constant} of the manifold.

Let us denote
$$\begin{array}{ll}
4\pi (X,Y)Z &:= g(Y,Z)X - g(X,Z)Y - 2g(JX,Y)JZ\\
[2mm]
&+ g(JY,Z)JX - g(JX,Z)JY,\\
[2mm]
8\Phi(X,Y)Z&:=g(Y,Z)(\eta(X)\xi-\eta(JX)J\xi)
-g(X,Z)(\eta(Y)\xi-\eta(JY)J\xi)\\
[2mm]
&+g(JY,Z)(\eta(X)J\xi+\eta(JX)\xi)-g(JX,Z)(\eta(Y)J\xi
+\eta(JY)\xi)\\
[2mm]
&-2g(JX,Y)(\eta(Z)J\xi+\eta(JZ)\xi)\\
[2mm]
&+(\eta(Y)\eta(Z)+
\eta(JY)\eta(JZ))X
-(\eta(X)\eta(Z)+\eta(JX)\eta(JZ))Y\\
[2mm]
&-(\eta(Y)\eta(JZ)-\eta(JY)\eta(Z))JX
+(\eta(X)\eta(JZ)-\eta(JX)\eta(Z))JY\\
[2mm]
&+2(\eta(X)\eta(JY)-\eta(JX)\eta(Y))JZ,
\quad X, Y, Z \in {\X}M.
\end{array}$$
In \cite{GM2} we have also proved that
\vskip 1mm
{\it If $(M,g,J) \; (\dim M = 2n \geq 6)$ is a Bochner-K\"ahler manifold
whose scalar distribution $D_{\tau}$ is a
$B_0$-distribution, then
$$R = a\pi + b\Phi, \quad b \neq 0,\leqno{(3.4)}$$
where $a, \, b$ are the following functions on $M$
$$a=\frac{\tau}{(n+1)(n+2)}+\frac{2{\bb}_0}{n+2},\quad
b=\frac{2\tau}{(n+1)(n+2)}-\frac{2n{\bb}_0}{n+2},\leqno{(3.5)}$$
and
$$b_0= \frac{2a-b}{2}= \, const. \leqno{(3.6)}$$}

In \cite{GM2} we studied three classes of Bochner-K\"ahler manifolds
whose scalar distribution is a $B_0$-distribution according to the function
$a+k^2$:
$$a+k^2>0,\quad a+k^2=0,\quad a+k^2<0.$$

Now we can prove a curvature characterization of K\"ahler manifolds
admitting a flat complex conformal connection.

\begin{thm} \label{T:3.1}
A K\"ahler manifold $(M,g,J) \; (\dim M = 2n \geq 6)$
admits a flat complex conformal connection if and only if
it is a Bochner-K\"ahler manifold whose scalar distribution
$D_{\tau}$ is a $B_0$-distribution with function $a+k^2=0$ and
geometric constants ${\B}=\frak b_0=0.$
\end{thm}
{\bf Proof.} Let $u$ be a $\mathcal{C}^{\infty}$ function on $M$,
such that the complex conformal connection $\mathcal{D}$, given
by (2.2) with $\omega = du\neq 0,\; P= grad \,u$, is flat.
Then (2.5) and (2.6) imply that the curvature tensor $R$ of
$M$ has the structure (1.1).
Comparing the tensor $Q$ from (1.1) and the tensor $L$ from (2.6)
we obtain
$$\rho(X,Y)=-2(n+2)\{\omega(X)\omega(Y)+\omega(JX)\omega(JY)+
\omega(P)g(X,Y)\},\; X, Y \in {\X}M \leqno(3.7)$$
and
$$\rho(X,P)=-2(n+2)\omega(P)\omega(X), \quad X\in {\X}M. \leqno(3.8)$$
After taking a trace in (3.7) we also get
$$\tau=-4(n+1)(n+2)\omega(P).\leqno (3.9)$$

Taking into account (2.5) we calculate from (3.7)
$$\begin{array}{ll}
(\nabla_X\,\rho)(Y,Z)=&2(n+2)\, \omega(P)\{2\omega(X)g(Y,Z)
+\omega(Y)g(X,Z)+\omega(Z)g(X,Y)\\
[3mm]
&+\omega(JY)g(X,JZ)+\omega(JZ)g(X,JY)\}, \quad X,Y,Z \in {\X}M.
\end{array}\leqno (3.10)$$
Comparing (3.1) and (3.10) in view of (3.9), we obtain
$$ \omega=-\frac{d\tau}{2\tau},\quad P=-\frac{grad\,\tau}{2\tau},
\quad \Vert d\tau\Vert^2=\frac{-{\tau}^3}{(n+1)(n+2)}.
\leqno (3.11)$$

The unit vector field $\displaystyle{\xi = \frac{grad \,
\tau}{\Vert d\tau \Vert}}$ because of (3.11) gets the form
$$\xi = 2\sqrt{\frac{(n+1)(n+2)}{-\tau}}\,P.$$

From (2.5) and (3.9) we obtain
$$\nabla_X\,\xi=-\frac{1}{2}\,\sqrt{\frac{-\tau}{(n+1)(n+2)}}
\{X-\eta(X)\xi-2\eta(JX)J\xi\}, \quad X \in {\X}M.
\leqno (3.12)$$

Now from (3.2) and (3.12) it follows that the scalar distribution
$D_{\tau}$ of the manifold is a $B_0$-distribution with functions
$$k=-\sqrt{\frac{-\tau}{(n+1)(n+2)}}\,, \quad
p^*=\frac{3}{2}\,\sqrt{\frac{-\tau}{(n+1)(n+2)}}\,.\leqno(3.13)$$
Then (2.6) and (3.11) give that the curvature tensor $R$ of the
manifold has the form
$$R=\frac{\tau}{(n+1)(n+2)}\,(\pi + 2\,\Phi)$$
and the functions $a$ and $b$ are
$$a=\frac{\tau}{(n+1)(n+2)}, \quad b=\frac{2\tau}{(n+1)(n+2)}.
\leqno(3.14)$$
From (3.13) and (3.14) we find $a+k^2=0$. The equalities (3.6) and (3.14)
imply that $\frak{b}_0=0$.

Taking into account (2.5), (3.11) and (3.7) we find
$$\Delta\,\tau=\frac{-\tau^2}{n+1}, \quad \Vert \rho \Vert^2
=\frac{(n+3)\tau^2}{2(n+1)^2}.\leqno (3.15)$$
Replacing $\Delta\,\tau$ and $\Vert \rho \Vert^2$ in (3.3) we obtain
${\B}=0$.
\vskip 2mm
For the inverse, let $(M,g,J)$ be a Bochner-K\"ahler manifold whose
scalar distribution is a $B_0$-distribution. Then it follows
\cite{GM2} that (3.2), (3.4) and (3.5) hold good. Under the condition
$\frak{b}_0=0$ we find that the functions $a$ and $b$ satisfy (3.14).

The condition $a+k^2=0$ implies that $k^2=-a=
\displaystyle{\frac{-\tau}{(n+1)(n+2)}}$\,,\, i.e. \,$\tau<0$. From
Theorem 3.5 in \cite{GM1} it follows that
$$\Vert d\tau \Vert = \xi(\tau)=\frac{(n+1)(n+2)}{2}\,\xi(b)=
\frac{(n+1)(n+2)}{2}\,kb>0,$$
which gives that the function $k$ is negative and
$$k=-\sqrt{\frac{-\tau}{(n+1)(n+2)}},\quad
\Vert d\tau\Vert^2=\frac{-{\tau}^3}{(n+1)(n+2)},\quad
p^* = \frac{3}{2}\,\sqrt{\frac{-\tau}{(n+1)(n+2)}}\,.$$
Then, from the equality (3.2) for any $X,Y \in {\X}M$ we have
$$(\nabla_X\,\eta)(Y)=-\frac{1}{2}\,\sqrt{\frac{-\tau}{(n+1)(n+2)}}
\{g(X,Y)-\eta(X)\eta(Y)+2\eta(JX)\eta(JY)\}.$$

Putting $2u:=-\ln {(-\tau)}$ and $\omega:=du=
-\displaystyle{\frac{d\tau}{2\tau}
=-\frac{\Vert d\tau \Vert}{2\tau}\,\eta=-\frac{k}{2}\,\eta}$
we prove that $(\nabla_X\,\omega)(Y)$ satisfies (2.5) and the
complex conformal connection (2.2) is flat.
\hfill {\bf QED}
\vskip 2mm
Let $(Q_0,g_0,\varphi,\tilde\xi_0, \tilde\eta_0)$ be an
$\alpha_0$-Sasakian space form \cite{JV} with constant
$\varphi$-holomorphic sectional curvatures $H_0$.
In \cite{GM2} we introduced warped product K\"ahler manifolds,
which are completely determined by the underlying
$\alpha_0$-Sasakian space form $Q_0$ of type
$H_0+3\alpha_0^2 \gtreqqless 0$
and the generating function $p(t), \; t \in I \subset \mathbb{R}$.

In order to obtain a local description of the K\"ahler manifolds
admitting a flat complex conformal connection we apply Theorem 6.1
in \cite {GM2}, which states:

{\it Any Bochner-K\"ahler manifold whose scalar distribution is a
$B_0$-distribution locally has the structure of a warped
product K\"ahler manifold with generating function $p(t)$
$($or $t(p))$ of type $1. - 13.$}

According to Theorem \ref{T:3.1} any K\"ahler manifold $(M,g,J),\;
(\dim M =2n \geqq 6)$ admitting a flat complex conformal connection
is locally a Bochner-K\"ahler manifold whose scalar distribution
is a $B_0$-distribution with function $a+k^2=0$ and constants
${\B}=\frak b_0=0$. In terms of \cite{GM2} the conditions
${\B}=\frak b_0=0$ are equivalent to the conditions
$\frak K= \frak b_0=0$.

Hence, $(M,g,J)$ is a warped product
Bochner-K\"ahler manifold whose underlying
$\alpha_0$-Sasakian space form is of type
$H_0+3\alpha_0^2 = 0$ with metric
$$g=p^2(t)\displaystyle{\left\{g_0+
\left(\frac{1}{\alpha_0}\frac{dp}{dt}-1\right)\,
\tilde\eta_0\otimes\tilde\eta_0\right\}+
\eta\otimes\eta},$$
generated by the function
$$p(t)=\frac{1}{\sqrt[3]{1-3\alpha_0(t-t_0)}}\,, \quad
t\in \left(-\infty,\,\frac{1+3\alpha_0t_0}{3\alpha_0}\right)$$
of type 9. \cite{GM2}.

This metric is not complete.

Especially in the case $\alpha_0=1$ the underlying manifold is a
Sasakian space form with $H_0=-3$. Sasakian space forms of type
$H_0=-3$ have been studied by Ogiue \cite {O} and Okumura \cite{Ok}.
A classification theorem for Sasakian space forms under the assumption
of completeness has been given by Tanno \cite{T}.

\end{document}